\documentclass[a4paper,11pt]{scrartcl}

\usepackage[T1]{fontenc}
\usepackage[latin1]{inputenc}
\usepackage[english]{babel}
\usepackage{setspace}
\usepackage{amsmath,amssymb,amsthm}
\usepackage{amsthm}
\usepackage{hyperref}
\usepackage{mathtools}
\usepackage{dsfont}
\usepackage{url}

\setkomafont{title}{\normalfont\bfseries}
\setkomafont{section}{\normalfont\bfseries}

\theoremstyle{plain}
\newtheorem{thm}{Theorem}[section]
\newtheorem{lem}[thm]{Lemma}
\theoremstyle{definition}
\newtheorem{defn}{Definition}[section]

\newcommand{\runinhead}[1]{\vspace{1mm}\noindent\textbf{#1}\hspace{0.5ex}}

\newcommand{\hnS}{\hspace{-1.25pt}}
\newcommand{\divergence}{\nabla\hnS\hnS\cdot\hnS}
\newcommand{\gradient}{\nabla\hnS}
\newcommand{\gradienth}{\nabla\hnS_h}
\newcommand{\boundary}[1]{\ensuremath{\partial\hnS #1}}
\newcommand{\restrict}[2]{{\ensuremath{\left. #1\right|_{#2}}}}
\newcommand{\komma}{\text{,}}
\newcommand{\punkt}{\text{.}}
\newcommand{\mydot}{\hnS\cdot\hnS}
\newcommand{\N}{\mathbb{N}}
\newcommand{\R}{\mathbb{R}}

\newcommand{\Pk}{\mathbb{P}}

\newcommand{\Vh}{V_h}
\newcommand{\Hdiv}{H_\text{div}}
\newcommand{\RTN}{RTN^l}

\newcommand{\param}{\mu}
\newcommand{\Params}{\mathcal{P}}
\newcommand{\paramFixed}{\overline{\mu}}
\newcommand{\eps}{\varepsilon}
\newcommand{\preeps}[1]{\prescript{\eps}{}{\hnS#1}}
\newcommand{\preparam}[1]{\prescript{}{\param}{#1}}
\newcommand{\preparamFixed}[1]{\prescript{}{\paramFixed}{#1}}
\newcommand{\preparamFixedBoth}[1]{\prescript{\eps}{\paramFixed}{#1}}
\newcommand{\preboth}[1]{\prescript{\eps}{\param}{#1}}
\newcommand{\prexi}[1]{\prescript{}{\xi}{#1}}

\newcommand{\alphaParam}{\prescript{}{\param, \paramFixed}{\alpha}}
\newcommand{\gammaParam}{\prescript{}{\param, \paramFixed}{\gamma}}

\newcommand{\Triangulation}{\mathcal{T}_H}
\newcommand{\triangulation}{\tau_h}
\newcommand{\Element}{T}
\newcommand{\element}{t}
\newcommand{\vertex}{\nu}

\newcommand{\Neighbour}{S}
\newcommand{\Faces}{\mathcal{F}_H}
\newcommand{\faces}{\mathcal{F}_h}
\newcommand{\face}{e}
\newcommand{\Face}{E}
\newcommand{\mean}[1]{\ensuremath{\left\{\hnS\hnS\left\{#1\right\}\hnS\hnS\right\}}}
\newcommand{\jump}[1]{\ensuremath{\left[\hnS\left[#1\right]\hnS\right]}}
\newcommand{\energynorm}[2]{\ensuremath{{\left|\hnS\left|\hnS\left|#1\right|\hnS\right|\hnS\right|}_{#2}}}
\newcommand{\norm}[2]{\ensuremath{{\left|\hnS\left|#1\right|\hnS\right|}_{#2}}}

\newcommand{\one}{\mathds{1}}
\DeclareMathOperator{\dx}{dx}

\title{A-posteriori error estimates for the localized reduced basis multi-scale method}
\author{Mario Ohlberger and Felix Schindler\thanks{Mario Ohlberger, Felix Schindler (formerly Albrecht), Applied Mathematics, University of M\"unster, Einsteinstr. 62, D-48149 M\"unster, \url{felix.schindler@wwu.de}}}
\date{}

\begin{document}

\maketitle

\begin{abstract}
We present a localized a-posteriori error estimate for the \underline{l}ocalized \underline{r}educed \underline{b}asis \underline{m}ulti-\underline{s}cale (LRBMS) method \cite{AHKO2012}.
The LRBMS is a combination of numerical multi-scale methods and model reduction using reduced basis methods to efficiently reduce the computational complexity of parametric multi-scale problems with respect to the multi-scale parameter $\eps$ and the online parameter $\param$ simultaneously.
We formulate the LRBMS based on a generalization of the SWIPDG discretization presented in \cite{ESV2010} on a coarse partition of the domain that allows for any suitable discretization on the fine triangulation inside each coarse grid element.
The estimator is based on the idea of a conforming reconstruction of the discrete diffusive flux, presented in \cite{ESV2010}, that can be computed using local information only.
It is offline/online decomposable and can thus be efficiently used in the context of model reduction.
\end{abstract}

\section{Introduction}
\label{fvca7-schindler-lrbms::section::introduction}

We are interested in efficient and reliable numerical approximations of parametric elliptic
multi-scale problems for given parameters $\param\in\Params \subset\mathbb{R}^p$, for $p \in \mathbb{N}$, i.e.
\begin{eqnarray}
  -\divergence(\preparam{\lambda} \preeps{\kappa} \mydot \gradient\preboth{p}) = f
    &&\text{in } \Omega\komma
  \label{fvca7-schindler-lrbms::equation::analytical_setting}
\end{eqnarray}
with homogeneous Dirichlet boundary conditions, where $\eps$ indicates the multi-scale nature of the quantities in prefix notation.
Equation (\ref{fvca7-schindler-lrbms::equation::analytical_setting}) arises e.g. in the context of two-phase flow in porous media, where it needs to be solved in every timestep for different $\param$ to obtain the global pressure $\preboth{p} :\Omega \to \R$  (see \cite[Sect. 1]{AHKO2012}).
A discretization of \eqref{fvca7-schindler-lrbms::equation::analytical_setting} usually consists in finding an approximation $\preboth{p_h} \in \Vh$ by a Galerkin projection onto a fine triangulation $\tau_h$ of $\Omega$ resolving the $\eps$ scale.
Two traditional approaches exist to reduce the computational complexity of the discrete problem: numerical multi-scale methods and model order reduction techniques.
Numerical multi-scale methods reduce the complexity of multi-scale problems with respect to $\eps$, while model order reduction techniques reduce the complexity of parametric problems with respect to $\mu$ (see \cite{Ohl2012} for an overview).
It is well known that solving parametric heterogeneous multi-scale problems accurately can be challenging and computationally costly, in particular for strongly varying scales and parameter ranges.
In general, numerical multi-scale methods capture the macroscopic behavior of the solution in a coarse approximation space, e.g., $V_H \subset V_h$, usually associated with a coarse triangulation $\Triangulation$ of $\Omega$, and recover the microscopic behavior of the solution by local fine-scale corrections.
Model order reduction using \underline{r}educed \underline{b}asis (RB) methods, on the other hand, is based on the idea to introduce a reduced space $V_{\text{red}} \subset V_h$, spanned by solutions of \eqref{fvca7-schindler-lrbms::equation::multiscale_dg_approximation} for a limited number of parameters $\mu$.
These training parameters are iteratively selected by an adaptive Greedy procedure (see \cite{AHKO2012} and the reference therein).
The idea of the recently presented \underline{l}ocalized \underline{r}educed \underline{b}asis \underline{m}ulti-\underline{s}cale (LRBMS) approach (see \cite{AHKO2012}) is to combine numerical multi-scale and RB methods and to generate a local reduced space $V_{\text{red}}^\Element \subset V_h^\Element$ for each coarse element of $\Triangulation$, given a tensor product type decomposition of the fine approximation space, $V_h = \oplus_{\Element\in\Triangulation}V_h^\Element$.
The coarse reduced space is then given as $V_{H,\text{red}} := \oplus_{\Element\in\Triangulation} V_{\text{red}}^\Element \subset V_h$, resulting in a multiplicative decomposition of the solution into $\preboth{p_{H,\text{red}}(x)}=\sum_{n=1}^{\dim(V_{H, \text{red}})}\preparam{p_n}(x)\preeps{\varphi_n}(x)$, where the RB functions $\preeps{\varphi_n}$ capture the microscopic behavior of the solution and the coefficient functions $\preparam{p_n}$ only vary on the coarse triangulation.

It is vital for an efficient and reliable use of RB as well as LRBMS methods to have access to an estimate on the model reduction  error.
Such an estimate is used to drive the adaptive Greedy basis generation during the offline phase of the computation and to ensure the quality of the reduced solution during the online phase.
It is usually given by a residual based estimator involving the stability constant and the residual in a dual norm.
It was shown in \cite{AHKO2012} that such an estimator can be successfully applied in the context of the LRBMS, but it was also pointed out that an estimator relying on global information might not be computationally feasible since too much work is required in the offline part of the computation.

The novelty of this contribution lies in a completely different approach to error estimation -- at least in the context of RB methods.
We make use of the ansatz of local error estimation presented in \cite{ESV2010} which measures the error by a conforming reconstruction of the physical quantities involved, specifically the diffusive flux $-\preparam{\lambda}\preeps{\kappa}\gradient \preboth{p}$.
This kind of local error estimation was proven to be very successful in the context of multi-scale problems and robust with respect to $\eps$.
We show in this work how we can transfer those ideas to the framework of the LRBMS to obtain an estimate of the error
$\preparam{\big|\hnS\big|\hnS\big| \preboth{p} - \preboth{p_{H,\text{red}}} \big|\hnS\big|\hnS\big|}$.
We would like to point out that we are able to estimate the error against the weak solution $\preboth{p}$ in a parameter dependent energy norm while traditional RB-approaches only allow to estimate the model reduction error in a parameter independent norm and only against the discrete solution.
In principal, this approach is able to turn the LRBMS method into a full multi-scale approximation scheme, while traditional RB methods can only be seen as a model reduction technique.
We would also like to point out that, to the best of our knowledge, this is the first work that makes use of local error information in the context of RB methods.

This work is organized as follows.
Section \ref{fvca7-schindler-lrbms::section::problem_formulation_discretization} introduces the notation and presents the overall setting, the discretization and the LRBMS framework.
We then carry out the error analysis for our multi-scale SWIPDG discretization in the parametric setting as well as the LRBMS method in Sec. \ref{fvca7-schindler-lrbms::section::error_analysis} and state our main result in Thm. \ref{fvca7-schindler-lrbms::theorem::locally_computable_abstract_energy_norm_estimate}.

\section{Problem formulation, discretization and model reduction}
\label{fvca7-schindler-lrbms::section::problem_formulation_discretization}

We consider linear elliptic problems of the form \eqref{fvca7-schindler-lrbms::equation::analytical_setting} in a bounded connected domain $\Omega \subset \R^d$, $d=2, 3$, with polygonal boundary $\boundary{\Omega}$ for a set of admissible parameters $\Params\subset\R^p$, $p\in\N$.

\runinhead{Triangulations} We require two nested partitions  of $\Omega$, a coarse one, $\Triangulation$, and a fine one, $\triangulation$.
Let $\triangulation$ be a simplicial triangulation of $\Omega$ with elements $\element \in \triangulation$.
In the context of multi-scale problems we call $\tau_h$ a \emph{fine triangulation} if it resolves all features of the quantities involved in \eqref{fvca7-schindler-lrbms::equation::analytical_setting}, specifically if $\preeps{\kappa}^\element := \restrict{\preeps{\kappa}}{\element}\in [L^\infty(\element)]^{d\times d}$ is constant for all $\element \in \triangulation$.
We only require the coarse elements $\Element\in\Triangulation$ to be shaped such that a local Poincar\'{e} inequality in $H^1(\Element)$ is fulfilled (see Thm. \ref{fvca7-schindler-lrbms::theorem::locally_computable_abstract_energy_norm_estimate}) and collect in $\triangulation^\Element \subset \triangulation$ the fine elements of $\triangulation$ that cover the coarse element $\Element$.
In addition, we collect all fine faces in $\faces$, all coarse faces in $\Faces$ and denote by $\Faces^\Element \subset \Faces$ the faces of a coarse element $\Element\in\Triangulation$ and by $\faces^\Face \subset \faces$ the fine faces that cover a coarse face $\Face\in\Faces$.

\runinhead{The continuous problem} We define the \emph{broken Sobolev space} $H^1(\triangulation) \subset L^2(\Omega)$ by $H^1(\triangulation) := \big\{ q \in L^2(\Omega) \;\big|\; q|_\element \in H^1(\element) \;\; \forall \element\in\triangulation \big\}$,
with $H^1_0(\Omega) \subset H^1(\Omega) \subset H^1(\triangulation)$, where $H^1$ denotes the usual Sobolev space of weakly differentiable functions and $H^1_0$ its elements which vanish on the boundary in the sense of traces.
In the same manner we denote the local broken Sobolev spaces $H^1(\triangulation^\Element)\subset L^2(\Element)$ for all $\Element\in\Triangulation$.
We also denote by $\gradienth:H^1(\triangulation)\to[L^2(\Omega)]^d$ the \emph{broken gradient operator} which is locally defined by $\restrict{(\gradienth q)}{\element}:=\gradient (\restrict{q}{\element})$ for all $\element\in\triangulation$ and $q\in H^1(\triangulation)$.
Given $f\in L^2(\Omega)$, $\preparam{\lambda} \in C^0(\Omega)$ strictly positive and $\preeps{\kappa} \in [L^\infty(\Omega)]^{d\times d}$ symmetric and uniformly positive definite, such that $\preparam{\lambda} \preeps{\kappa}\in [L^\infty(\Omega)]^{d\times d}$ is bounded from below (away from 0) and above for all $\param\in\Params$, we define the parameter dependent bilinear form $\preboth{b}: H^1(\triangulation) \times H^1(\triangulation) \to \R$ and the linear form $l: H^1(\triangulation) \to \R$ by $\preboth{b}(p,q) := \sum_{\Element\in\Triangulation} \preboth{b}^T(p,q)$ and $l(q) := \sum_{\Element\in\Triangulation} l^\Element(q)$, respectively, and their local counterparts $\preboth{b}^\Element := H^1(\triangulation^\Element) \times H^1(\triangulation^\Element) \to \R$ and $l^\Element: H^1(\triangulation^\Element) \to \R$ for all $\Element\in\Triangulation$ and $\param\in\Params$ by
\begin{equation}
  \preboth{b}^\Element(p,q) := \int\limits_\Element(\preparam{\lambda}\preeps{\kappa}\mydot \gradienth p)\mydot\gradienth q \dx
    \quad\quad\text{ and }\quad\quad
  l^\Element(q):= \int\limits_\Element f q\dx.
  \nonumber
\end{equation}

\begin{defn}[Weak solution]
  \label{fvca7-schindler-lrbms::definition::weak_solution}
  Given a parameter $\param\in\Params$ we define the \emph{weak solution} of \eqref{fvca7-schindler-lrbms::equation::analytical_setting} by $\preboth{p}\in H^1_0(\Omega)$, such that
  \begin{equation}
    \preboth{b}(\preboth{p},q) = l(q) \quad\quad\text{for all } q \in H^1_0(\Omega)\punkt
    \label{fvca7-schindler-lrbms::equation::weak_solution}
  \end{equation}
\end{defn}

Note that, since $\preboth{b}$ is continuous and coercive for all $\param\in\Params$ (due to the assumptions on $\preparam{\lambda}\preeps{\kappa}$) and since $l$ is bounded, there exists a unique solution of \eqref{fvca7-schindler-lrbms::equation::weak_solution} due to the Lax-Milgram Theorem.

\runinhead{A note on parameters} In addition to the assumptions we posed on $\preparam{\lambda}$ above we also demand it to be \emph{affinely decomposable} with respect to $\param\in\Params$, i.e. there exist $\varXi\geq 1$ strictly positive \emph{coefficients} $\prexi{\theta}:\Params\to\R$ for $0\leq\xi\leq\varXi -1$ and $\prescript{}{\varXi}{\theta}\in\{0,1\}$ and $\varXi + 1$ nonparametric \emph{components} $\prexi{\lambda}\in C^0(\Omega)$, such that $\preparam{\lambda} = \sum_{\xi=0}^{\varXi} \prexi{\theta}(\param) \prexi{\lambda}$.
We can then compare $\lambda$ for two parameters $\param, \paramFixed \in\Params$ by $\prescript{}{\param, \paramFixed}{\alpha}\preparamFixed{\lambda} \;\leq\; \preparam{\lambda} \;\leq\; \prescript{}{\param, \paramFixed}{\gamma} \preparamFixed{\lambda}$, where $\prescript{}{\param, \paramFixed}{\alpha} :=\min_{\xi = 0}^{\varXi - 1} {\prexi{\theta}(\param)} {\prexi{\theta}({\paramFixed})}^{-1}$ and $\prescript{}{\param, \paramFixed}{\gamma} := \max_{\xi = 0}^{\varXi - 1} {\prexi{\theta}(\param)}{\prexi{\theta}({\paramFixed})}^{-1}$
denote the positive equivalence constants.
This assumption on the data function $\preparam{\lambda}$ is a common assumption in the context of RB methods and covers a wide range of physical problems.
If $\preparam{\lambda}$ does not exhibit such a decomposition one can replace $\preparam{\lambda}$ by an arbitrary close approximation using Empirical Interpolation techniques (see \cite{AHKO2012} and the references therein) which does not impact our analysis.
All quantities that linearly depend on $\preparam{\lambda}$ inherit the above affine decomposition in a straightforward way.
Since we would like to estimate the error in a problem dependent norm we also need the notion of a \emph{parameter dependent energy norm} $\preboth{\energynorm{\cdot}{}}: H^1(\triangulation)\to\R$ for any $\param\in\Params$, defined by $\preboth{\energynorm{q}{}} := \big(\sum_{\Element\in\Triangulation}\preboth{\energynorm{q}{\Element}}^2\big)^{1/2}$ with $\preboth{\energynorm{q}{T}} := \big(\preboth{b}^\Element(q,q)\big)^{1/2}$,
for all $\Element\in\Triangulation$.
Note that $\preboth{\energynorm{\cdot}{}}$ is a norm only on $H^1_0(\Omega)$.
We can compare these norms for any two parameters $\param, \paramFixed\in\Params$ using the above decomposition of $\preparam{\lambda}$:
\begin{equation}
  \sqrt{\prescript{}{\param, \paramFixed}{\alpha}} \;\; \preparamFixedBoth{\energynorm{\cdot}{}}
    \;\;\;\leq\;\;\; \preboth{\energynorm{\cdot}{}}
    \;\;\;\leq\;\;\; \sqrt{\prescript{}{\param, \paramFixed}{\gamma}} \;\; \preparamFixedBoth{\energynorm{\cdot}{}}
  \label{fvca7-schindler-lrbms::equation::norm_equivalence}
\end{equation}
We denote by $0 < \preboth{c}_\element \leq \preboth{C}_\element$ the smallest and largest eigenvalue of $\preparam{\lambda}^\element \preeps{\kappa}^\element$ and additionally define $0 < \preeps{c}^\element := \min_{\param\in\Params} \preboth{c}^\element$, $\preeps{c}^\element < \preeps{C}^\element := \max_{\param\in\Params} \preboth{C}^\element$ for all $\element\in\triangulation$.
From here on we denote an a-priori chosen parameter by $\hat{\param} \in \Params$ while $\paramFixed\in\Params$ denotes an arbitrary parameter and $\param \in\Params$ denotes the parameter during the online phase of the simulation.

\runinhead{The generalized SWIPDG discretization} We discretize \eqref{fvca7-schindler-lrbms::equation::weak_solution} by allowing for a suitable discretization of at least first order inside each coarse element $\Element\in\Triangulation$ and by coupling those with a \underline{s}ymmetric \underline{w}eighted \underline{i}nterior \underline{p}enalty (SWIP) \underline{d}iscontinuous \underline{G}alerkin (DG) discretization along the coarse faces of $\Triangulation$.
We give a very brief definition of the SWIPDG bilinear form, see \cite[Sect. 2.3]{ESV2010} and the references therein for a detailed discussion and the definition of $\jump{\cdot}_\face$, $\mean{\cdot}_\omega$ and $\gamma_\face$.
For any two-valued function $q\in H^1(\triangulation)$, we define its multi-scale jump $\jump{q}_\Face$ and its multi-scale average $\preeps{\mean{q}_\Face}$ for all coarse faces $\Face\in\Faces$ locally by $\restrict{\jump{q}_{\Face}}{\face} := \jump{q}_\face$ and $\restrict{\preeps{\mean{q}_\Face}}{\face} := \mean{q}_\omega$ for all $\face\in\faces^\Face$.
In addition we define the multi-scale penalty parameter $\preboth{\sigma_\Face}$ locally by $\restrict{\preboth{\sigma_\Face}}{\face} := \preparam{\lambda} \gamma_\face$ for all fine faces $\face\in\faces^\Face$ on all coarse faces $\Face\in\Faces$.
With these definitions at hand we define the \emph{multi-scale SWIPDG} bilinear form $\preboth{b_h}: H^1(\triangulation)\times H^1(\triangulation)\to\R$ by
\begin{equation}
  \preboth{b_h}(p,q)
    := \sum\limits_{\Element\in\Triangulation} \preboth{b_h^\Element}(p,q)
    + \sum\limits_{\Face\in\Faces} \preboth{b_h^\Face}(p,q)\komma
  \label{fvca7-schindler-lrbms::equation::multiscale_swipdg_bilinear_forms}
\end{equation}
with the coupling bilinear forms $\preboth{b_h}^\Face: H^1(\triangulation^\Element)\times H^1(\triangulation^\Neighbour)\to\R$ given by
\begin{eqnarray}
  \preboth{b_h^\Face}(p,q)
    :=\int\limits_{\Face}
      &-&\preeps{\mean{(\preparam{\lambda} \preeps{\kappa}\mydot\gradienth q)\mydot n_\Face}_\Face} \jump{p}_\Face
  \nonumber\\
      &-&\Big(\preeps{\mean{(\preparam{\lambda} \preeps{\kappa}\mydot\gradienth p)\mydot n_\Face}_\Face}
      -\preboth{\sigma_\Face}\jump{p}_\Face\Big)\jump{q}_\Face
      \dx
  \nonumber
\end{eqnarray}
for all $\Face = \boundary\Element \cap \boundary\Neighbour\in\Faces$.
To complete the definition of the discretization we only demand the local bilinear forms $\preboth{b_h^\Element}$ to be an approximation of $\preboth{b^\Element}$ (with trivial Neumann boundary values) and the local discrete ansatz spaces $\Vh^{k,\Element}$ to be locally polynomial of order $k \geq 1$, i.e. $\restrict{q}{\element} \in \Pk_k(\element)$ for all $\element\in\triangulation^\Element$ and $q\in\Vh^{k, \Element}$ on all $\Element\in\Triangulation$.
We then define the multi-scale DG approximation space as $\Vh^k(\Triangulation) := \big\{ q\in H^1(\triangulation) \big| \restrict{q}{\Element} \in \Vh^{k,\Element} \;\;\forall\Element\in\Triangulation \big\} \subset H^1(\triangulation)$ for $k \geq 1$.
\begin{defn}[Multi-scale DG approximation]
\label{fvca7-schindler-lrbms::definition::multiscale_dg_approximation}
  Given a parameter $\param\in\Params$ we define the \emph{multi-scale DG approximation} of \eqref{fvca7-schindler-lrbms::equation::weak_solution} by $\preboth{p_h}\in \Vh^1(\Triangulation)$, such that
  \begin{equation}
    \preboth{b_h}(\preboth{p_h}, q_h) = l(q_h) \quad\quad\text{for all } q_h \in \Vh^1(\Triangulation)\punkt
    \label{fvca7-schindler-lrbms::equation::multiscale_dg_approximation}
  \end{equation}
\end{defn}
The bilinear form $\preboth{b_h}$ is continuous and coercive if the penalty parameter is chosen large enough and if the sum of the local bilinear forms is continuous and coercive.
If those are chosen accordingly the discrete problem \eqref{fvca7-schindler-lrbms::equation::multiscale_dg_approximation} thus has a unique solution.
Possible choices for the local bilinear forms $\preboth{b_h^\Element}$ and the local approximation spaces $\Vh^{k,\Element}$ include continuous Finite Elements and variants of the IPDG and the SWIPDG discretizations.

\runinhead{The localized reduced basis multi-scale method} Since the global (in a spatial sense) model reduction ansatz of classical RB methods does not always fit in the context of multi-scale problems, the LRBMS introduced in \cite{AHKO2012} takes a more localized approach to model reduction.
We refer to \cite{AHKO2012} for a detailed definition of the LRBMS and only state what is needed for the error analysis here.
The main idea of the LRBMS is to restrict solutions of \eqref{fvca7-schindler-lrbms::equation::multiscale_dg_approximation} for some $\param$ to the elements of the coarse triangulation and to form local reduced spaces $V_{\text{red}}^\Element \subset \Vh^{k, \Element}$ by a local compression of those solution snapshots.
Given these local reduced spaces we define the broken reduced space by $V_{H,\text{red}} := \oplus_{\Element\in\Triangulation} V_{\text{red}}^\Element \subset \Vh^k(\Triangulation)$.
The LRBMS approximation is then given by a standard Galerkin projection of \eqref{fvca7-schindler-lrbms::equation::multiscale_dg_approximation}.
\begin{defn}[LRBMS approximation]
\label{fvca7-schindler-lrbms::definition::LRBMS_approximation}
  Given a parameter $\param\in\Params$ we define the \emph{LRBMS approximation} of \eqref{fvca7-schindler-lrbms::equation::weak_solution} by $\preboth{p_{H,\text{red}}}\in V_{H,\text{red}}$, such that
  \begin{equation}
    \preboth{b_h}(\preboth{p_{H,\text{red}}}, q_H) = l(q_H) \quad\quad\text{for all } q_H \in V_{H,\text{red}}\punkt
    \label{fvca7-schindler-lrbms::equation::LRBMS_approximation}
  \end{equation}
\end{defn}

\section{Error analysis}
\label{fvca7-schindler-lrbms::section::error_analysis}

Our error analysis is a generalization of the ansatz presented in \cite{ESV2010} to provide an estimator for our multi-scale DG approximation solving \eqref{fvca7-schindler-lrbms::equation::multiscale_dg_approximation} as well as for our LRBMS approximation solving \eqref{fvca7-schindler-lrbms::equation::LRBMS_approximation}.
We transfer the idea of a conforming reconstruction of the nonconforming discrete diffusive flux $-\preparam{\lambda}\preeps{\kappa}\gradienth p_h$ to our setting.
Our error analysis shares some similarities with the general multi-scale ansatz presented in \cite{PVWW2013}, which is stated for a wide range of discretizations but for a different coupling strategy.

We obtain the mild requirement for the local approximation spaces that the constant function $\one$ is present, which is obvious for traditional discretizations and can be easily achieved for the LRBMS approximation by incorporating the DG basis with respect to $\Triangulation$.
The estimates are fully offline/online decomposable and can thus be used for efficient model reduction in the context of the LRBMS.

We begin by stating an \emph{abstract energy norm estimate} (see \cite[Lemma 4.1]{ESV2010}) that splits the difference between the weak solution $\preboth{p}\in H^1_0(\Omega)$ solving \eqref{fvca7-schindler-lrbms::equation::weak_solution} and any function $p_h\in H^1(\triangulation)$ into two contributions.
This abstract estimate does not depend on our discretization and thus leaves the choice of $s$ and $u$ open.
Note that we formulate the following Lemma with separate parameters for the energy norm and the weak solution.
The price we have to pay for this flexibility are the additional constants involving $\alphaParam$ and $\gammaParam$, that vanish if $\paramFixed$ and $\param$ coincide.
\begin{lem}[Abstract energy norm estimate]
\label{fvca7-schindler-lrbms::lemma::abstract_energy_norm_estimate}
  Given any $\param,\paramFixed\in\Params$ let $\preboth{p}\in H^1_0(\Omega)$ be the weak solution solving \eqref{fvca7-schindler-lrbms::equation::weak_solution} and let $p_h\in H^1(\triangulation)$ be arbitrary.
  Then
  \begin{eqnarray}
    \preparamFixedBoth{\energynorm{\preboth{p} - p_h}{}}
      &\leq& \tfrac{1}{\sqrt{\prescript{}{\param, \paramFixed}{\alpha}}} \Big(
      \inf\limits_{s\in H^1_0(\Omega)} \sqrt{\gammaParam} \preparamFixedBoth{}{\energynorm{p_h - s}{}}
    \nonumber\\
    &+& \inf\limits_{u\in \Hdiv(\Omega)}
      \Big\{
        \sup\limits_{\substack{\varphi\in H^1_0(\Omega)\\\preboth{\energynorm{\varphi}{}}=1}}
          \big\{
            \big(
              f - \divergence u, \varphi
            \big)_{L^2}
            -\big(
              \preparam{\lambda}\preeps{\kappa}\mydot\gradienth p_h + u, \gradient\varphi
            \big)_{L^2}
          \big\}
      \Big\}\Big)
    \nonumber\\
    &\leq& \tfrac{\sqrt{\gammaParam}}{\sqrt{\alphaParam}} \;\;2\;\; \preparamFixedBoth{\energynorm{\preboth{p} - p_h}{}}\punkt
    \nonumber
  \end{eqnarray}
\end{lem}
The above Lemma is proven by applying the norm equivalence \eqref{fvca7-schindler-lrbms::equation::norm_equivalence}, following the arguments in the proof of \cite[Lemma 4.1]{ESV2010} and applying the norm equivalence again.

The next Thm. states the main localization result and gives an indication on how to proceed with the choice of $u$:
it allows us to localize the estimate of the above Lemma, if $u \in \Hdiv(\Omega) := \big\{ v \in [L^2(\Omega)]^{d\times d} \;\big|\; \divergence v \in L^2(\Omega)\big\}$ fulfills a \emph{local conservation property}.

\begin{thm}[Locally computable abstract energy norm estimate]
\label{fvca7-schindler-lrbms::theorem::locally_computable_abstract_energy_norm_estimate}
\mbox{}
  Let $\preboth{p}\in H^1_0(\Omega)$ be the weak solution of \eqref{fvca7-schindler-lrbms::equation::weak_solution}, let $s\in H^1_0(\Omega)$ and $p_h\in H^1(\triangulation)$ be arbitrary, let $u\in\Hdiv(\Omega)$ fulfill the \emph{local conservation property} $(\divergence{u}, \mathds{1})_\Element = (f, \mathds{1})_\Element$ and let $C_P^\Element > 0$ denote the constant from the Poincar\'{e} inequality $\norm{\varphi - \Pi_0^\Element\varphi}{L^2,\Element}^2 \leq C_P^\Element h_\Element^2 \norm{\gradient\varphi}{L^2,\Element}^2$ for all $\varphi\in H^1(\Element)$ on all $\Element\in\Triangulation$, where $\Pi_l^\omega$ denotes the $L^2$-orthogonal projection onto $\Pk_l(\omega)$ for $l \in \N$ and $\omega \subseteq \Omega$.
  It then holds that
  \begin{eqnarray}
    \preparamFixedBoth{\energynorm{\preboth{p} - p_h}{}}
      \;\;&\leq&\;\; \eta[s, u]\komma
      \quad\quad\text{with the \emph{global estimator} $\eta[s,u]$ defined as}
    \nonumber\\
    \eta[s,u]
      &:=& \tfrac{\sqrt{\gammaParam}}{\sqrt{\alphaParam}}
        \Big(\sum_{\Element\in\Triangulation} \eta_\text{nc}^\Element[s]^2\Big)^{1/2}
    \nonumber\\
      &+&\tfrac{1}{\sqrt{\alphaParam}}\Big(\sum_{\Element\in\Triangulation} \eta_\text{r}^\Element[u]^2\Big)^{1/2}
      +\tfrac{\max\big(\sqrt{\prescript{}{\mu, \hat{\mu}}{\gamma}}, \sqrt{\prescript{}{\mu, \hat{\mu}}{\alpha}}^{-1}\big)}{\sqrt{\alphaParam \prescript{}{\param, \hat{\param}}{\alpha}}} \Big(\sum_{\Element\in\Triangulation} \eta_\text{df}^\Element[u]^2\Big)^{1/2}
    \nonumber
  \end{eqnarray}
  and the local \emph{nonconformity estimator} given by $\eta_\text{nc}^\Element[s] := \preparamFixedBoth{\energynorm{p_h - s}{\Element}}$, the local \emph{residual estimator} given by $\eta_\text{r}^\Element[u] := ({C_P^\Element}/{\preeps{c}^\Element})^{1/2} h_\Element \norm{f - \divergence u}{L^2,\Element}$ and the local \emph{diffusive flux estimator} given by $\eta_\text{df}^\Element[u] := \norm{(\prescript{}{\hat{\param}}{\lambda}\preeps{\kappa})^{1/2}\gradienth p_h + (\prescript{}{\hat{\param}}{\lambda}\preeps{\kappa})^{-1/2} u}{L^2, \Element}$ for all coarse elements $\Element\in\Triangulation$, where $\preeps{c^\Element} := (\max_{\element\in\triangulation^\Element} 1/{\preeps{c^\element}})^{-1}$.
\end{thm}
The above Thm. is proven by loosely following the proof of \cite[Thm. 3.1]{ESV2010}, i.e. by starting from Lem. \ref{fvca7-schindler-lrbms::lemma::abstract_energy_norm_estimate}, localizing with respect to $\mathcal{T}_H$, using the local conservation property and the norm equivalence \eqref{fvca7-schindler-lrbms::equation::norm_equivalence}.

What is left now in order to turn the abstract estimate of Thm. \ref{fvca7-schindler-lrbms::theorem::locally_computable_abstract_energy_norm_estimate} into a fully computable one is to specify $s$ and $u$.
We will do so in the following paragraphs.

\runinhead{Oswald interpolation} Given any nonconforming approximation $p_h \in \Vh^k(\Triangulation) \not\subset H^1_0(\Omega)$ we will choose $s \in H^1_0(\Omega)$ as a conforming reconstruction of $p_h$
by the \emph{Oswald Interpolation operator} $I_\text{os} : \Vh^1(\Triangulation) \to \Vh^1(\Triangulation) \cap H^1_0(\Omega)$ which we define by prescribing its values on the Lagrange nodes of the triangulation (see \cite[Sect. 2.5]{ESV2010} and the references therein):
we define $I_\text{os}[p_h](\vertex) := p_h^\element(\nu)$ inside any $\element\in\triangulation$ and
\begin{equation}
  I_\text{os}[p_h](\vertex) := \tfrac{1}{|\triangulation^v|} \sum_{\element \in \triangulation^\vertex} p_h^\element(\vertex)
    \quad\text{for all inner nodes of } \triangulation \text{ and}\quad I_\text{os}[p_h](\vertex) := 0
  \nonumber
\end{equation}
for all boundary nodes of $\triangulation$, where $\triangulation^v \subset \triangulation$ denotes the set of all simplices of the fine triangulation which share $\vertex$ as a node.

\runinhead{Diffusive flux reconstruction} As mentioned above we will reconstruct a conforming diffusive flux approximation $u_h \in \Hdiv(\Omega)$ of the nonconforming discrete diffusive flux $-\preparam{\lambda}\preeps{\kappa}\gradienth p_h \not\in \Hdiv(\Omega)$ in a conforming discrete subspace $\RTN_h(\triangulation) \subset \Hdiv(\Omega)$, namely the \emph{Raviart-Thomas-N\'{e}d\'{e}lec} space of vector functions (see \cite{ESV2010} and the references therein), defined for $k - 1 \leq l \leq k$ by
\begin{equation}
  \RTN_h(\triangulation)
    :=\big\{
        v \in \Hdiv(\Omega)
      \big|
        \restrict{v}{\element} \in \RTN(\element) := [\Pk_l(\element)]^d + \vec{x} \Pk_l(\element)
        \quad\forall\element\in\triangulation
      \big\}\punkt
  \nonumber
\end{equation}
See \cite[Sect. 2.4]{ESV2010} and the references therein  for a detailed discussion of the role of the polynomial degree $l$, the properties of elements of $\RTN_h(\triangulation)$ and the origin of the use of diffusive flux reconstructions in the context of error estimation in general.
Now, given any $p_h \in H^1(\triangulation)$ and any $\param\in\Params$ we define the diffusive flux reconstruction $\preboth{u_h}[p_h]\in\RTN_h(\triangulation)$ locally by demanding
\begin{equation}
  \Big(
      \preboth{u_h}[p_h]\mydot n_\Face, q
    \Big)_{L^2, \Face}
    =\Big(
      -\preeps{\mean{(\preparam{\lambda}\preeps{\kappa}\mydot\gradienth p_h)\mydot n_\Face}_\Face}
        + \sigma_\Face \jump{p_h}_\Face, q
    \Big)_{L^2, \Face}
  \nonumber
\end{equation}
for all $q\in\Pk_l(\face)$ for all $\face\in\faces^\Face$ and all $\Face\in\Faces^\Element$ and by
\begin{equation}
  \Big(
      \preboth{u_h}[p_h], \gradienth q
    \Big)_{L^2, \Element}
    =- \preboth{b_h^\Element}(p_h, q)
    +\sum_{\Face\in\Faces^\Element}
      \Big(
        \preeps{\omega_\Face^+} (\preparam{\lambda}^\Element\preeps{\kappa}^\Element\mydot \gradienth q)\mydot n_\Face, \jump{p_h}_\Face
      \Big)_{L^2, \Face}
  \nonumber
\end{equation}
for all $q \in \Vh^{k, \Element}$ such that $\restrict{\gradient q}{\element} \in [\Pk_{l-1}(t)]^d$ for all $\element\in\triangulation^\Element$ and all $\Element\in\Triangulation$.
The next Lemma shows that this reconstruction of the diffusive flux is sensible for a multi-scale approximation as well as an LRBMS approximation, since the reconstructions of both fulfill the requirements of Thm. \ref{fvca7-schindler-lrbms::theorem::locally_computable_abstract_energy_norm_estimate}.

\begin{lem}[Local conservativity]
  \label{fvca7-schindler-lrbms::lemma::local_conservativity}
  Let $\preboth{p_*}\in H^1(\triangulation)$ either denote a multi-scale DG approximation $\preboth{p_h} \in \Vh^1(\Triangulation)$ given by \eqref{fvca7-schindler-lrbms::equation::multiscale_dg_approximation} or an LRBMS approximation $\preboth{p_{H,\text{red}}} \in V_{H,\text{red}}$ given by \eqref{fvca7-schindler-lrbms::equation::LRBMS_approximation}. Let $\preboth{u_h}[\preboth{p_*}]\in\RTN_h(\triangulation)$ denote its diffusive flux reconstruction and let $\one \in V^{*, \Element}$, where $V^{*, \Element}$ either denotes the local approximation space $\Vh^{1,\Element}$ or the local reduced space $V_{\text{red}}^\Element$, for all $\Element\in\Triangulation$.
  Then $\preboth{u_h}[\preboth{p_*}]$ fulfills the local conservation property of Thm. \ref{fvca7-schindler-lrbms::theorem::locally_computable_abstract_energy_norm_estimate}.
\end{lem}
The above Lemma is proven by applying the ideas of \cite[Lemma 2.1]{ESV2010} to our setting while accounting for $\mathcal{T}_H$, i.e. by using the local conservation property, the definition of the discrete bilinear form and the fact, that $\mathds{1} \in V^{*, {T}}$.
At this points some remarks are in order.
If we drop the parameter dependency and set $\Triangulation = \triangulation$, we obtain the discretization proposed in \cite{ESV2010} and the estimators of Thm. \ref{fvca7-schindler-lrbms::theorem::locally_computable_abstract_energy_norm_estimate} and \cite[Thm. 3.1]{ESV2010} coincide.
The estimators defined in Thm. \ref{fvca7-schindler-lrbms::theorem::locally_computable_abstract_energy_norm_estimate} can be efficiently offline/online decomposed, even if we choose $\paramFixed = \param$.
A more elaborate work containing the proofs and the efficiency of the estimator (using standard arguments) is in preparation.

We finally obtain a fully computable and fully specified estimate by combining the definition of the Oswald interpolant and the diffusive flux reconstruction with Thm. \ref{fvca7-schindler-lrbms::theorem::locally_computable_abstract_energy_norm_estimate} for both our multi-scale DG discretization and the LRBMS method.

{\small
\bibliographystyle{plaindin}
\bibliography{estimates-lrbms-2014-ohlberger-schindler}
}
\end{document}